\documentclass[12pt]{article}

\usepackage{amsmath}
\usepackage{amsthm}
\usepackage{amssymb}
\usepackage{cite}	            
\usepackage{url}
\usepackage[multiple]{footmisc} 
\usepackage{graphicx}
\usepackage{epstopdf}
\usepackage{chngcntr}           
\usepackage{enumitem}
\usepackage{hhline}
\usepackage{array}
\usepackage{hyperref}
\usepackage{color}

\hypersetup{colorlinks=false}

\usepackage[pagewise]{lineno} 

\if@twoside
\else 
\fi

\numberwithin{equation}{section}
\counterwithin{figure}{section}
\counterwithin{table}{section}

\theoremstyle{plain}
\newtheorem{theorem}{Theorem}[section]

\newtheorem{algorithm}  [theorem]{Algorithm}

\theoremstyle{definition}

\theoremstyle{remark}
\newtheorem*{remark}{Remark}

\newcommand{\dis}{\displaystyle}
\newcommand{\norm}[1]{\left\|#1\right\|}

\newcommand{\spr}[1]{\left\langle\,#1\,\right\rangle}
\newcommand{\kl}[1]{\left(#1\right)}

\newcommand{\hsp}{\hspace*{5mm}}
\newcommand{\hhsp}{\hspace*{1cm}}
\newcommand{\hhhsp}{\hspace*{15mm}}

\newcommand{\R}{\mathbb{R}}

\newcommand{\ddt}{\frac{\partial}{\partial t}}

\newcommand{\mm}{\text{mm}}

\newcommand{\s}{\text{s}}

\newcommand{\citeV}{\cite{Voss_Dyke_Tabelow_Schiff_Ballon_2017}}

\newcommand{\citeVH}{\cite{Voss_Dyke_Tabelow_Schiff_Ballon_2017,Hubmer_Neubauer_Ramlau_Voss_2018}}

\newcommand{\Oijk}{\Omega_{i,j,k}}
\newcommand{\rb}{\bar{r}}
\newcommand{\LtO}{{L^2(\Omega)}}

\title{A Conjugate-Gradient Approach to the Parameter Estimation Problem of Magnetic Resonance Advection Imaging}
\author{
Simon Hubmer\footnote{Johann Radon Institute Linz, Altenbergerstra{\ss}e 69, A-4040 Linz, Austria, (simon.hubmer@ricam.oeaw.ac.at), Corresponding author} ,
Andreas Neubauer\footnote{Johannes Kepler University Linz, Institute of Industrial Mathematics, Altenbergerstra{\ss}e 69, A-4040 Linz, Austria (neubauer@indmath.uni-linz.ac.at)} ,
Ronny Ramlau\footnote{Johannes Kepler University Linz, Institute of Industrial Mathematics, Altenbergerstra{\ss}e 69, A-4040 Linz, Austria, (ronny.ramlau@jku.at)} \footnote{Johann Radon Institute Linz, Altenbergerstra{\ss}e 69, A-4040 Linz, Austria, (ronny.ramlau@ricam.oeaw.ac.at)}\,,
Henning U. Voss\footnote{Weill Cornell Medical College, Department of  Radiology, 516 E 72nd Street
New York, NY 10021 (hev2006@med.cornell.edu)}}

\begin{document}

\maketitle

\begin{abstract}
We consider the inverse problem of estimating the spatially varying pulse wave velocity in blood vessels in the brain from dynamic MRI data, as it appears in the recently proposed imaging technique of Magnetic Resonance Advection Imaging (MRAI). The problem is formulated as a linear operator equation with a noisy operator and solved using a conjugate gradient type approach. Numerical examples on experimental data show the usefulness and advantages of the developed algorithm in comparison to previously proposed methods.

\end{abstract}



\section{Introduction}

With every beat of the human heart a pulse wave is created which travels through the blood vessels in the body somewhat like the waves produced by a rock thrown into the water.  The speed of the this wave, termed pulse wave velocity (PWV), is closely related to material properties of the blood vessel it travels through \cite{Li_2004,Moens_1877,Korteweg_1878}; a high PWV indicates a stiffer blood vessel than a lower PWV, given the same vessel wall thickness and vessel diameter. Arterial pulse wave velocity (PWV) is the gold standard for aortic stiffness assessment in cardiovascular disease \cite{Laurent_Cockcroft_VanBortel_Boutouyrie_Giannattasio_Hayoz_Pannier_Vlachopoulos_Wilkinson_Struijker-Boudier_2006, Mancia_2007, Wentland_Grist_Wieben_2014} and provides normative values for healthy and increased arterial stiffness 
\cite{McEniery_Yasmin_Hall_Qasem_Wilkinson_Cockcroft_ACCT_Investigators_2005, Boutouyrie_Vermersch_Laurent_Briet_2008, Mattace_Raso_et_al_2010}. 
Arterial stiffness is related to arterial compliance \cite{Elter_2001, Bramwell_1922,Vulliemoz_Stergiopulos_Meuli_2002}, the ratio of blood volume change to blood pressure change. Arterial compliance absorbs the impact of the pulse pressure waves and enables steady blood flow throughout the whole cardiac cycle \cite{Guyton_Hall_2006, Vrselja_Brkic_Mrdenovic_Radic_Curic_2014}. Pulse wave velocity and arterial compliance are reliable prognostic markers for cardiovascular morbidity and mortality in adult populations such as the elderly, subjects with diabetes, arteriosclerosis, coronary heart disease, and hypertension \cite{Kozakova_Morizzo_Guarino_Federico_Miccoli_Giannattasio_Palombo_2015, Macgowan_Stoops_Zhou_Cahill_Sled_2015, Bogren_Mohiaddin_Klipstein_Firmin_Underwood_Rees_Longmore_1989, Heintz_Gillessen_Walkenhorst_Dahl_Dorr_Krebs_Hanrath_Sieberth_1993, Mohiaddin_Firmin_Longmore_1993, Honda_Yano_Matsuoka_Hamada_Hiwada_1994, Rogers_Hu_Coast_Vido_Kramer_Pyeritz_Reichek_2001, Meaume_Benetos_Henry_Rudnichi_Safar_2001, McEniery_Yasmin_Hall_Qasem_Wilkinson_Cockcroft_ACCT_Investigators_2005, Boutouyrie_Vermersch_Laurent_Briet_2008, Mattace_Raso_et_al_2010}.

It would be desirable if measurements of the PWV could not only be performed for the major arteries in the body but also in the brain using Magnetic Resonance Imaging (MRI). Brain MRI already provides high-resolution angiograms of arteries deep inside the brain. Adding the pulsatile component of blood flow could allow for the measurement of the cerebral PWV along cerebral arteries. The long-term goal would be to provide a novel clinical biomarker to assess cerebrovascular integrity, namely the cerebrovascular PWV.
As a first attempt to non-invasively measure the PWV using MRI, Voss et al. \citeV recently introduced the technique of Magnetic Resonance Advection Imaging (MRAI). They propose to measure the pulsatile component of blood flow
from dynamic echo-planar imaging (EPI) data, such as acquired in functional and resting-state functional MRI experiments \cite{Stehling_Turner_Mansfield_1991,Ogawa_Lee_Kay_Tank_1990}. The fundamental idea is that the PWV and the associated MRI data can be connected via the advection equation
	\begin{equation}\label{eq_transport}
		\ddt \rho(x,y,z,t) + v(x,y,z) \cdot \nabla \rho(x,y,z,t) = 0 \,,
	\end{equation}
where $\rho(x,y,z,t)$ denotes the time-dependent MRI data and $v(x,y,z)$ the spatially varying PWV, which is assumed to be divergence-free. Equation \eqref{eq_transport} was derived from physical considerations and multiple regression was used to estimate $v(x,y,z)$. While leading to promising first results, for example a correct identification of the pulse wave direction following the main cerebral arteries, it was determined that this approach is not suitable to quantitatively estimate the PWV. The main reason is insufficient accuracy of the EPI data and limited temporal resolution. 

Hence, in \cite{Hubmer_Neubauer_Ramlau_Voss_2018} the authors proposed a stable and efficient algorithm for estimating $v(x,y,z)$ via \eqref{eq_transport} based on a formulation of this parameter estimation problem within the framework of Inverse Problems \cite{Engl_Hanke_Neubauer_1996}. The resulting nonlinear inverse problem was solved using a fast gradient-based iterative regularization method, which yielded stable reconstructions of the PWV in a reasonable computational time for the given spatio-temporal resolution of the considered MRI data sets. Unfortunately, as with the regression approach, only qualitative estimates of the PWV could be obtained, since the spatio-temporal resolution of the MRI data set was too low to yield quantitative results. However, it is expected that the resolution of MRI data will increase in the years to come, hopefully allowing also quantitative estimates of the PWV. While the algorithm proposed in \cite{Hubmer_Neubauer_Ramlau_Voss_2018} is, in principle, able to deal with data sets of higher resolution, it may, in practice, not be able to handle the resulting large amount of data in a reasonable time. Thus, there is the need for a faster and robust algorithm for estimating the PWV, which scales well with respect to the size of the involved data sets.

In this paper, we present an algorithm satisfying the above requirements, which is based on a reformulation of the problem into a linear system of equations with both a noisy operator and a noisy right-hand side. This system is then solved by the method of Conjugate Gradients for the Normal Equations (CGNE), where the required calculations can be written down explicitly in terms of simple calculations without having to assemble or store the system matrix, which allows for effective implementation. The remaining part of this paper is structured as follows: in Section~\ref{sect_mathematical_model}, we derive the exact mathematical model and reformulate the problem into a system of linear equations, based on a suitable space-time discretization of the problem. In Section~\ref{section_CGNE}, we consider the application of CGNE to the problem, deriving an explicit algorithm to carry out the required computations, and which avoids the assembly of the system matrix. In Section~\ref{sect_numerical_results}, we present numerical results of the application of our algorithm to experimental data, and finally, in Section~\ref{sect_conclusion}, we give some short conclusions. 

Since the aim of this paper is mainly the design and presentation of an effective algorithm for solving the parameter estimation problem of MRAI, the interested reader is referred to \citeVH for a detailed medical and physical background on MRAI.

\section{Mathematical Model and Discretization}\label{sect_mathematical_model}

In this section, we consider a suitable mathematical model and numerical discretization for the parameter estimation problem of MRAI, leading to a linear inverse problem with both a noisy operator and a noisy right-hand side. 

As originally proposed in \citeV, the physics-based model of MRAI implied that the MRI signal $\rho$ is a conserved quantity advected (transported) by the velocity vector field $v$. In its general form, the advection (transport, continuity) equation is given by 
	\begin{equation}\label{advection}
		\rho_t(x,y,z,t)+\nabla\cdot(\rho(x,y,z,t)\,v(x,y,z))=0\,,\qquad (x,y,z)\in \Omega\,,t\in[0,T]\,.
	\end{equation}
Under the assumption of a divergence-free vector field $v$, which is reasonable for an incompressible flow, this leads to model equation \eqref{eq_transport} used in \citeVH. However, since this assumption does not necessarily hold in practice, we directly work with equation \eqref{advection} here.

A typical MRI data set does not provide continuous information over the whole area of the brain but rather consists of averaged values given on a voxel grid. Hence, we assume (for simplicity) that the domain $\Omega$ corresponding to the observed region of the brain is a cuboid. Furthermore, since the voxels are typically uniform, we set
	\begin{equation*}
 		x_i:=x_0+i\Delta x\,,\quad y_j:=y_0+j\Delta y\,,\quad z_k:=z_0+k\Delta z\,,
	\end{equation*}
and define the voxels $\Oijk:=[x_{i-1},x_i]\times[y_{j-1},y_j]\times[z_{k-1},z_k]$. Then it follows from integrating \eqref{advection} over the sets $\Oijk$ that
	\begin{equation}\label{advv}
 		\frac{\partial}{\partial t}\int_{\Oijk}\rho(x,y,z,t)\,d(x,y,z)+
 		\int_{\partial\Oijk}\rho(x,y,z,t)\,v(x,y,z)\cdot n(x,y,z)\,dS=0\,,
	\end{equation}
where $n(x,y,z)$ denotes the unit outward normal and $1\le i\le I$, $1\le j\le J$, $1\le k\le K$.

Since neither $\rho$ nor $v$ are known everywhere, all integrals are approximated using some quadrature rules. If one uses midpoint rules, then the values at the boundaries of $\Oijk$ have to be approximated via interpolation. Therefore, we suggest the following rules that are exact for bilinear and trilinear functions, respectively:
	\begin{align*}
 		\int_{x_{i-1}}^{x_i}\int_{y_{j-1}}^{y_j}f(x,y)\,d(x,y) \approx {} & \frac{\Delta x
 		\Delta y}{4}\Big(f(x_{i-1},y_{j-1})+f(x_i,y_{j-1})
 		\\
 		& \hspace{1.5cm} + \, f(x_{i-1},y_j)+f(x_i,y_j)\Big),
 		\\
 		\int_{x_{i-1}}^{x_i}\int_{y_{j-1}}^{y_j}\int_{z_{k-1}}^{z_k}f(x,y,z)\,d(x,y,z)
 		\approx {} & \frac{\Delta x\Delta y\Delta z}{8}\Big(f(x_{i-1},y_{j-1},z_{k-1})+
 		f(x_i,y_{j-1},z_{k-1})\\
 		& \hspace{2cm}+\,f(x_{i-1},y_j,z_{k-1})+f(x_i,y_j,z_{k-1})\\
 		& \hspace{2cm}+\,f(x_{i-1},y_{j-1},z_k)+f(x_i,y_{j-1},z_k)\\
 		& \hspace{2cm}+\,f(x_{i-1},y_j,z_k)+f(x_i,y_j,z_k)\Big).
	\end{align*}
Thus, we obtain the semi-discrete system
	\begin{equation}\label{advsd}
 		\frac{\partial}{\partial t}D_{i,j,k}(t)=\frac{2}{\Delta x}A_{i,j,k}(t)+
 		\frac{2}{\Delta y}B_{i,j,k}(t)+\frac{2}{\Delta z}C_{i,j,k}(t)\,,
	\end{equation}
where $1\le i\le I$, $1\le j\le J$, $1\le k\le K$,
	\begingroup
	\allowdisplaybreaks
	\begin{align*}
		D_{i,j,k}(t) := {} & \rho_{i-1,j-1,k-1}(t)+\rho_{i,j-1,k-1}(t)+\rho_{i-1,j,k-1}(t)
 		+\rho_{i,j,k-1}(t)\\
 		& +\,\rho_{i-1,j-1,k}(t)+\rho_{i,j-1,k}(t)+\rho_{i-1,j,k}(t)+\rho_{i,j,k}(t)\,
 		,\\
		A_{i,j,k}(t) := {} & \Big(\rho_{i-1,j-1,k-1}(t)\,v_{1,i-1,j-1,k-1}+
		\rho_{i-1,j,k-1}(t)\,v_{1,i-1,j,k-1}\\
 		&\;\;+\,\rho_{i-1,j-1,k}(t)\,v_{1,i-1,j-1,k}+\rho_{i-1,j,k}(t)\,v_{1,i-1,j,k}
 		\Big)\\
 		&-\,\Big(\rho_{i,j-1,k-1}(t)v_{1,i,j-1,k-1}+\rho_{i,j,k-1}(t)v_{1,i,j,k-1}\\
 		&\quad\;\,+\,\rho_{i,j-1,k}(t)v_{1,i,j-1,k}+\rho_{i,j,k}(t)v_{1,i,j,k}\Big),\\
 		B_{i,j,k}(t) := {} & \Big(\rho_{i-1,j-1,k-1}(t)\,v_{2,i-1,j-1,k-1}+
 		\rho_{i,j-1,k-1}(t)\,v_{2,i,j-1,k-1}\\
 		&\;\;+\,\rho_{i-1,j-1,k}(t)\,v_{2,i-1,j-1,k}+\rho_{i,j-1,k}(t)\,v_{2,i,j-1,k}
 		\Big)\\
 		&-\,\Big(\rho_{i-1,j,k-1}(t)v_{2,i-1,j,k-1}+\rho_{i,j,k-1}(t)v_{2,i,j,k-1}\\
 		&\quad\;\,+\,\rho_{i-1,j,k}(t)v_{2,i-1,j,k}+\rho_{i,j,k}(t)v_{2,i,j,k}\Big),\\ 
 		C_{i,j,k}(t) := {}& \Big(\rho_{i-1,j-1,k-1}(t)\,v_{3,i-1,j-1,k-1}+
 		\rho_{i,j-1,k-1}(t)\,v_{3,i,j-1,k-1}\\
 		&\;\;+\,\rho_{i-1,j,k-1}(t)\,v_{3,i-1,j,k-1}+\rho_{i,j,k-1}(t)\,v_{3,i,j,k-1}
 		\Big)\\
 		&-\,\Big(\rho_{i-1,j-1,k}(t)v_{3,i-1,j-1,k}+\rho_{i,j-1,k}(t)v_{3,i,j-1,k}\\
 		&\quad\;\,+\,\rho_{i-1,j,k}(t)v_{3,i-1,j,k}+\rho_{i,j,k}(t)v_{3,i,j,k}\Big),
	\end{align*}
	\endgroup
and
	\begin{equation*}
		\rho_{i,j,k}(t):=\rho(x_i,y_j,z_k,t)\quad \text{and} \quad v_{m,i,j,k}:=  v_m(x_i,y_j,z_k)\,, 
	\end{equation*}
$0\le i\le I$, $0\le j\le J$, $0\le k\le K$, and $m=1,2,3$.

In a final step, the derivative with respect to $t$ is approximated by a
differential quotient. If one wants to calculate $\rho$ for given $v$, it is
more appropriate to use a forward differential quotient, since it is much faster
(but not as stable).

However, we are interested in calculating $v$ from measurements of $\rho$.
Therefore, we can as well use a backwards differential quotient. Unfortunately,
measurements of $\rho$ do not exist at all time steps for all grid points, due to what is called the \emph{slice-time} acquisition problem. The issue is that MRI data are usually gathered slice by slice, meaning that for each time step the data $\rho_{i,j,k}$ are only measured for a fixed $k$ depending on the time-step. Different scanning set-ups allow for different orders of the slices on which the data are acquired. As in \cite{Hubmer_Neubauer_Ramlau_Voss_2018}, we here focus (without loss of generality) on the case of ascending slice-time acquisition, i.e., measurements are available at points
	\begin{equation}\label{pt}
 		(x_i,y_j,z_k,t_{k,l})\,,\qquad 0\le i\le I,\,0\le j\le J,\,0\le k\le K,\,
 		0\le l\le L,
	\end{equation}
where
	\begin{equation*}
		t_{k,l}:=(k+(K+1)l)\Delta t\,.
	\end{equation*}  
This means that in each time step only one $z$-layer can be measured and after a
full cycle it restarts again. The measurements are abbreviated by $\rho_{i,j,k,l} := \rho(x_i,y_j,z_k,t_{k,l})$.

To reduce the size of the linear system obtained from \eqref{advsd}, we use the
following fully discrete system
	\begin{equation}\label{advd}
 		\frac{D_{i,j,k}(t_{k,l})-D_{i,j,k}(t_{k,l-1})}{(K+1)\Delta t}=
 		\frac{2}{\Delta x}A_{i,j,k}(t_{k,l})+\frac{2}{\Delta y}B_{i,j,k}(t_{k,l})+
 		\frac{2}{\Delta z}C_{i,j,k}(t_{k,l})\,,
	\end{equation}
where $1\le i\le I$, $1\le j\le J$, $1\le k\le K$, $1\le l<L$. Since measurements for $\rho$ are not available at all points, we have to use (linear) interpolation, i.e.,
	\begin{align*}
 		\rho_{i,j,k}(t_{k,l}) &=\rho_{i,j,k,l}\,,\hspace{7.2cm}0\le l\le L\,,\\
 		\rho_{i,j,k-1}(t_{k,l}) &= \rho_{i,j,k-1,l}\kl{1-\frac{1}{K+1}}+
 		\rho_{i,j,k-1,l+1}\frac{1}{K+1}\,,\qquad 0\le l<L\,.
	\end{align*}
If we also allow $l=L$ in \eqref{advd}, then $\rho_{i,j,k-1}(t_{k,L})$ has to be
approximated via extrapolation. Note that other slice-time acquisition procedures can be treated in a similar way using a suitably different interpolation scheme.

For our further calculations we assume that $\Delta x=\Delta y=\Delta z:=\Delta$.
Then the system \eqref{advd} may be written as the following discrete linear
equation
	\begin{equation}\label{lin}
		 Tv=b\,,	
	\end{equation}
where $v\in X:=\R^{3(I+1)(J+1)(K+1)}$ and $b\in Y:=\R^{I \cdot J \cdot K \cdot(L-1)}$ is defined by
	\begin{equation*}
		b_{i,j,k,l}:=\frac{(D_{i,j,k}(t_{k,l})-D_{i,j,k}(t_{k,l-1}))\Delta}{2(K+1)
		 \Delta t}\,.
	\end{equation*} 
According to \eqref{advd}
	\begin{equation*}
		(Tv)_{i,j,k,l}:=A_{i,j,k}(t_{k,l})+B_{i,j,k}(t_{k,l})+C_{i,j,k}(t_{k,l})\,.
	\end{equation*} 
One can solve equation \eqref{lin} using the CGNE method, as described in the next section.

\begin{remark}
Note that there are different possibilities to discretize the advection equation \eqref{advection}. For example, higher order schemes can often improve the quality of the solution to the forward problem, which, due to the influence of noise in the data, unfortunately does not necessarily translate to the inverse problem.

The advantage of the proposed discretization scheme is two-fold. On the one hand, it is well adapted to the structure and properties of the measured data. On the other hand, and more importantly, it ultimately allows to derive the explicit form of CGNE presented in Algorithms~3.1 and 3.2 below without actually having to compute the matrices $T$, $T^*$ or $T^*T$ explicitly. Extending the derivation to other (higher order) schemes might certainly be possible, but appears to be impracticable, due to the involved high computational effort. As reasonable results can also be obtained with our scheme, an extension in this direction is out of the scope of this work.
\end{remark}

\section{Solving the Discrete Inverse Problem using CGNE}\label{section_CGNE}

In this section, we consider CGNE for solving \eqref{lin}, which is based on the observation that a least squares solution of \eqref{lin} is given as a solution of the normal equations
	\begin{equation*}\label{lin_norm}
		T^* T v = T^* b \,,
	\end{equation*}
where $T^*$ denotes the adjoint of $T$. Since $T^*T$ is positive semi-definite, the method of conjugate gradients (CG) can be applied for solving \eqref{lin_norm}, which leads to CGNE.

This method has been extensively studied, especially in the framework of inverse problems with noisy data (see for example \cite{Engl_Hanke_Neubauer_1996} and the references therein). In particular, it is known that CGNE gives rise to a convergent regularization method if combined with a suitable stopping rule, and that the iterates converge to the minimum-norm solution.

For the CGNE method we need inner products in $X$ and $Y$. In $Y$ we choose the standard Euclidean inner product
	\begin{equation*}
		\spr{b,c}:=\sum_{i=1}^I\sum_{j=1}^J\sum_{k=1}^K\sum_{l=1}^{L-1}b_{i,j,k,l} c_{i,j,k,l}\,.
	\end{equation*} 
Since the functions $(v_1,v_2,v_3)$ are assumed to be differentiable, especially $v_1$ with respect to $x$, $v_2$ with respect to $y$, and $v_3$ with respect to $z$, we suggest the following inner product in $X$:
	\begin{align*}
		 \spr{v,w} &:= \sum_{i=0}^I\sum_{i=0}^J\sum_{k=0}^K(v_{1,i,j,k}w_{1,i,j,k}+
		 v_{2,i,j,k}w_{2,i,j,k}+v_{3,i,j,k}w_{3,i,j,k})\\
		 & +\,\frac{1}{\Delta^2}\sum_{i=1}^I\sum_{i=0}^J\sum_{k=0}^K(v_{1,i,j,k}-
		 v_{1,i-1,j,k})(w_{1,i,j,k}-w_{1,i-1,j,k})\\
		 & +\,\frac{1}{\Delta^2}\sum_{i=0}^I\sum_{i=1}^J\sum_{k=0}^K(v_{2,i,j,k}-
		 v_{2,i,j-1,k})(w_{2,i,j,k}-w_{2,i,j-1,k})\\
		 & +\,\frac{1}{\Delta^2}\sum_{i=0}^I\sum_{i=0}^J\sum_{k=1}^K(v_{3,i,j,k}-
		 v_{3,i,j,k-1})(w_{3,i,j,k}-w_{3,i,j,k-1}) \,, 
	\end{align*}
which is a discretized version of the inner product
	\begin{equation*}
		\spr{v,w}_\LtO \,+  \,\spr{\frac{\partial v_1}{d x} , \frac{\partial w_1}{d x}  }_\LtO 
		+ \, \spr{\frac{\partial v_2}{d y} , \frac{\partial w_2}{d y}  }_\LtO
		+ \,\spr{\frac{\partial v_3}{d z} , \frac{\partial w_3}{d z}  }_\LtO \,.
	\end{equation*}
We also need the adjoint of $T$ with respect to these inner products, i.e.,
	\begin{equation*}
		 \spr{Tv,d}_Y=\spr{v,T^*d}_X\,. 
	\end{equation*}
Note that
	\begin{align*}
		\lefteqn{\sum_{i=0}^I\sum_{i=0}^J\sum_{k=0}^Kv_{1,i,j,k}w_{1,i,j,k}+
		 \frac{1}{\Delta^2}\sum_{i=1}^I\sum_{i=0}^J\sum_{k=0}^K(v_{1,i,j,k}-
		 v_{1,i-1,j,k})(w_{1,i,j,k}-w_{1,i-1,j,k})}\\
		 &= \frac{1}{\Delta^2}\sum_{i=0}^J\sum_{k=0}^K (v_{1,0,j,k},\ldots,v_{1,I,j,k})
		 \kl{\begin{array}{rrrrr} a&-1&0&0&0\\-1&b&-1&0&0\\\ddots&\ddots&\ddots&\ddots&
		 \ddots\\ 0&0&-1&b&-1\\0&0&0&-1&a \end{array}}\kl{\begin{array}{c}w_{1,0,j,k}\\w_{1,1,j,k}\\
		 \vdots\\ w_{1,I-1,j,k}\\w_{1,I,j,k}\end{array}}	
	\end{align*}
with $a:=\Delta^2+1$ and $b:=a+1$. A similar formula holds for the parts concerning \linebreak $v_{2,\cdot,\cdot,\cdot}$ and $v_{3,\cdot,\cdot,\cdot}$.

For the calculation of the adjoint $T^*d$ we need the following numbers
	\begin{equation*}
		r_{-1}:=1\,,\qquad r_0:=a\,,\qquad r_i:=b\,r_{i-1}-r_{i-2}\,,\quad i=1,
		 \ldots,\max\{I,J,K\}\,,	
	\end{equation*} 
	\begin{equation*}
		\rb_I:=r_I-r_{I-1}\,,\qquad\rb_J:=r_J-r_{J-1}\,,\qquad\rb_K:=r_K-r_{K-1}\,.
	\end{equation*} 
Using the above equality, $T^*d$ may be calculated componentwise as follows:

\begin{algorithm}\label{tsd} Calculation of $T^*d$
\begin{itemize}\itemsep 0pt
	\item Set $\kappa:=0$
	\item Calculate: $c_{\alpha,\beta,\gamma,i,j,k}:=\dis\sum_{l=1}^{L-1}
 	d_{i+\alpha,j+\beta,k+\gamma,l}\rho_{i,j,k}(t_{k+\gamma,l})\,,\qquad\alpha,
 	\beta,\gamma\in\{0,1\}$,\\
 	\hspace*{2cm}$i=1-\alpha,\ldots,I-\alpha\,,\quad j=1-\beta,\ldots,J-\beta\,,\quad
 	k=1-\gamma,\ldots,K-\gamma$.\\
 	\hspace*{2cm}For all other indices $i,j,k$: $c_{\alpha,\beta,\gamma,i,j,k}:=0$.
	\item For $j=0,\ldots,J$ and $k=0,\ldots,K$ do:

	\hsp For $i=0,\ldots,I$ do:\\
	\hhsp $e_i:=c_{1,1,1,i,j,k}+c_{1,0,1,i,j,k}+c_{1,1,0,i,j,k}+c_{1,0,0,i,j,k}$\\
	\hspace*{2cm}$-\,c_{0,1,1,i,j,k}-c_{0,0,1,i,j,k}-c_{0,1,0,i,j,k}-c_{0,0,0,i,j,k}$\\
	\hhsp If ($i=0)$ \{\\
	\hhhsp $w_{1,i,j,k}:=e_i$\\
	\hhsp\} else \{\\
	\hhhsp $w_{1,i,j,k}:=e_i\,r_{i-1}+w_{1,i-1,j,k}$\\
	\hhsp\}

	\hhsp $w_{1,I,j,k}:=w_{1,I,j,k}/\rb_I\qquad\kappa:=\kappa+e_I*w_{1,I,j,k}$

	\hhsp For $i=I-1,\ldots,0$ do:\\
	\hhhsp $w_{1,i,j,k}:=(w_{1,i,j,k}+r_{i-1}w_{1,i+1,j,k})/r_i\qquad\kappa:=\kappa+
	e_i*w_{1,i,j,k}$
	\item For $i=0,\ldots,I$ and $k=0,\ldots,K$ do:

	\hsp For $j=0,\ldots,J$ do:\\
	\hhsp $e_j:=c_{1,1,1,i,j,k}+c_{0,1,1,i,j,k}+c_{1,1,0,i,j,k}+c_{0,1,0,i,j,k}$\\
	\hspace*{2cm}$-\,c_{1,0,1,i,j,k}-c_{0,0,1,i,j,k}-c_{1,0,0,i,j,k}-c_{0,0,0,i,j,k}$\\
	\hhsp If ($j=0)$ \{\\
	\hhhsp $w_{2,i,j,k}:=e_j$\\
	\hhsp\} else \{\\
	\hhhsp $w_{2,i,j,k}:=e_j\,r_{j-1}+w_{2,i,j-1,k}$\\
	\hhsp\}

	\hhsp $w_{2,i,J,k}:=w_{2,i,J,k}/\rb_J\qquad\kappa:=\kappa+e_J*w_{2,i,J,k}$

	\hhsp For $j=J-1,\ldots,0$ do:\\
	\hhhsp $w_{2,i,j,k}:=(w_{2,i,j,k}+r_{j-1}w_{2,i,j+1,k})/r_j\qquad\kappa:=\kappa+
 	e_j*w_{2,i,j,k}$
	\item For $i=0,\ldots,I$ and $j=0,\ldots,J$ do:

	\hsp For $k=0,\ldots,K$ do:\\
	\hhsp $e_k:=c_{1,1,1,i,j,k}+c_{1,0,1,i,j,k}+c_{0,1,1,i,j,k}+c_{0,0,1,i,j,k}$\\
	\hspace*{2cm}$-\,c_{1,1,0,i,j,k}-c_{1,0,0,i,j,k}-c_{0,1,0,i,j,k}-c_{0,0,0,i,j,k}$\\
	\hhsp If ($k=0)$ \{\\
	\hhhsp $w_{3,i,j,k}:=e_k$\\
	\hhsp\} else \{\\
	\hhhsp $w_{3,i,j,k}:=e_k\,r_{k-1}+w_{3,i,j,k-1}$\\
	\hhsp\}

	\hhsp $w_{3,i,j,K}:=w_{3,i,j,K}/\rb_K\qquad\kappa:=\kappa+e_K*w_{3,i,j,K}$

	\hhsp For $k=K-1,\ldots,0$ do:\\
	\hhhsp $w_{3,i,j,k}:=(w_{3,i,j,k}+r_{k-1}w_{3,i,j,k+1})/r_k\qquad\kappa:=\kappa+
 	e_k*w_{3,i,j,k}$
	\item Then $(T^*d)_{m,i,j,k}:=\Delta^2w_{m,i,j,k}$ and $\spr{T^*d,T^*d}_X:=
 	\Delta^2\kappa$.
\end{itemize}
\end{algorithm}

The CGNE method in algorithmic form is now given as follows:

\begin{algorithm}\label{cgne}{\bf(CGNE)}
\begin{itemize}\itemsep 0pt
\item Calculate the numbers $h:=\Delta^2$, $a:=h+1$, $b:=a+1$, and $r_i$,
 $\rb_I$, $\rb_J$, $\rb_K$ as above.
\item Set $v:=0$, $it:=0$, and choose $itmax$.
\item While $it<itmax)\,\{$

\hsp If (it = 0) \{

\hhsp $d:=b$\\
\hhsp Calculate $w$ and $\kappa$ as in Algorithm \ref{tsd}\\
\hhsp $\gamma:=\kappa$ and $p:=w$

\hsp\} else \{

\hhsp $d:=d-\alpha q$\\
\hhsp Calculate $w$ and $\kappa$ as in Algorithm \ref{tsd}\\
\hhsp $\beta:=\kappa/\gamma$, $\gamma:=\kappa$, and $p:=w+\beta p$

\hsp\}

\hsp Calculate $q:=Tp$\\
\hsp $\alpha:=\gamma/\spr{q,q}_Y$\\
\hsp $v:=v+\alpha p$\\
\hsp $it$++

\}
\end{itemize}
\end{algorithm}

\begin{remark}
Note that in Algorithm~\ref{cgne} the matrices $T$ and $T^*$ are not needed explicitly, which makes its application very time-efficient.
\end{remark}

\begin{remark}
In Algorithm~\ref{cgne}, the iteration is stopped after \emph{itmax} steps in order to avoid instabilities due to data noise, which are unavoidable due to the ill-posed nature of the problem. Thus, the stopped iteration only yields an approximation of the minimum-norm solution of the normal equation. Among other factors, the quality of this approximation depends on the stopping rule used to determine \emph{itmax}.

Unfortunately, many common choices are not suitable in our case for two reasons: firstly, no reliable estimate of the noise in the data is known, and secondly, it is not clear whether they are suitable in the case of a noisy operator. A partial remedy for this problem might be so-called \emph{noise-free} stopping rules, which, however, have not yet been thoroughly analysed for the case of noisy operators with unknown noise levels.

Thus, for our numerical problems described below, we have chosen \emph{itmax} via a comparison between the behaviour of the residual $\norm{Tv - b}$ and the image quality.
\end{remark}

\section{Numerical Results}\label{sect_numerical_results}

In this section, we consider the application of our CGNE method to experimental data. 

In order to be able to compare our methods to the existing techniques, we use the same MRI data sets as in \cite{Hubmer_Neubauer_Ramlau_Voss_2018}, namely the MRI scans of subjects $2$ and $16$ of a publicly available data set obtained on a $7.0$ T MRI scanner \cite{Hanke_Baumgartner_Ibe_Kaule_Pollmann_Speck_Zinke_Stadler_2014}. As in \cite{Hubmer_Neubauer_Ramlau_Voss_2018} we use the first $20$ seconds of the second $15$ min segments of the data. For these data sets we have that $\Delta = 1.4 \, \mm$ and $(K + 1)\Delta t = 2 \, \s$, which means that once every $2$ s the scanner has completed a single scan of the brain. 

As already mentioned in the last remark of the previous section, we have used a monitoring of the residual $\norm{Tv - b}$ and the image quality during the iteration to determine a suitable stopping index \emph{itmax}. The choice of $\text{\emph{itmax}} = 10$ was found to give the best results for our data sets, especially since the behaviour of the residual becomes unstable for higher numbers of iterations.

\begin{figure}[ht!]
\centering
\includegraphics[trim = 80mm 40mm 60mm 35mm, clip, width=\textwidth]{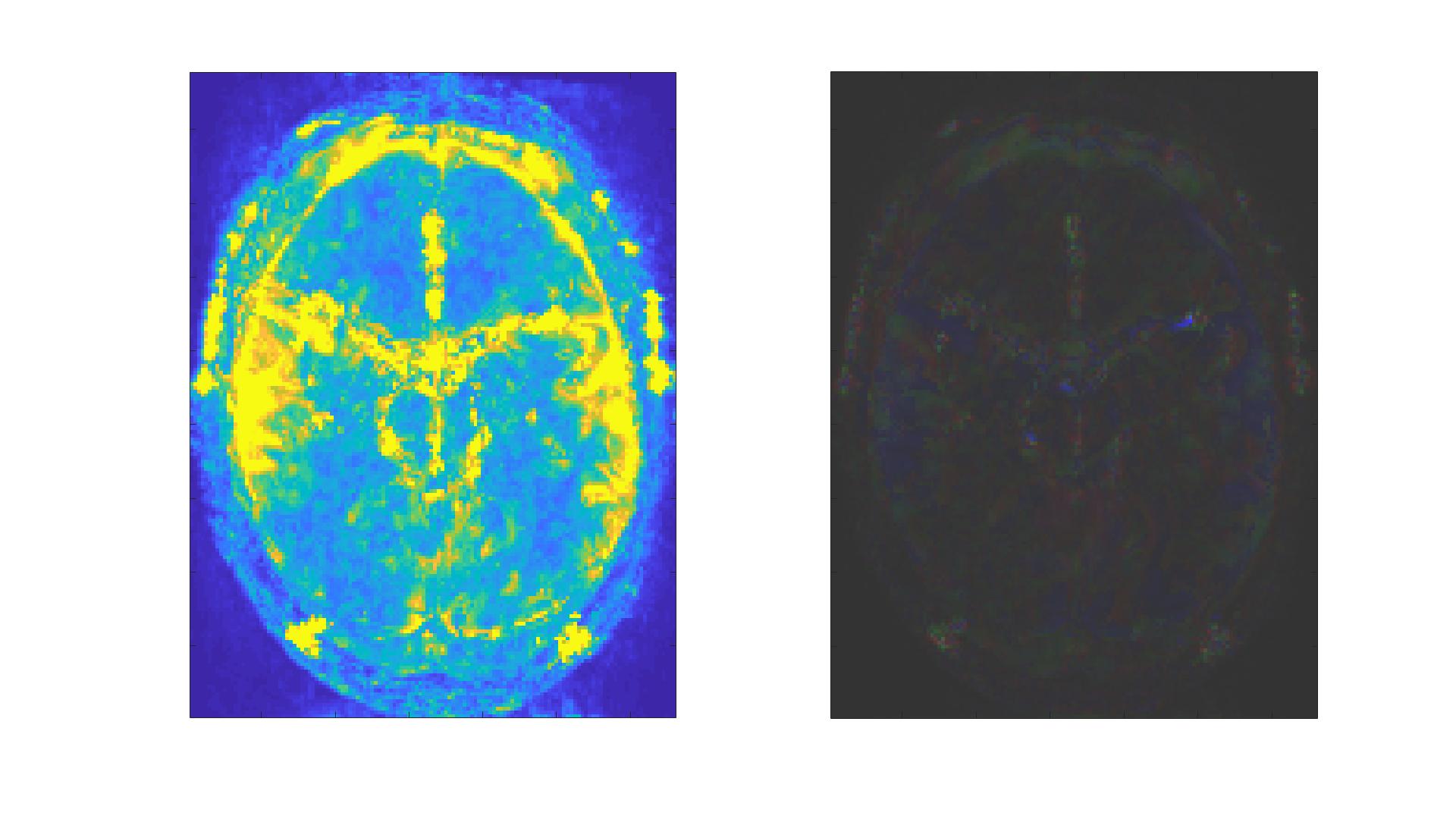}
\includegraphics[trim = 80mm 40mm 60mm 35mm, clip, width=\textwidth]{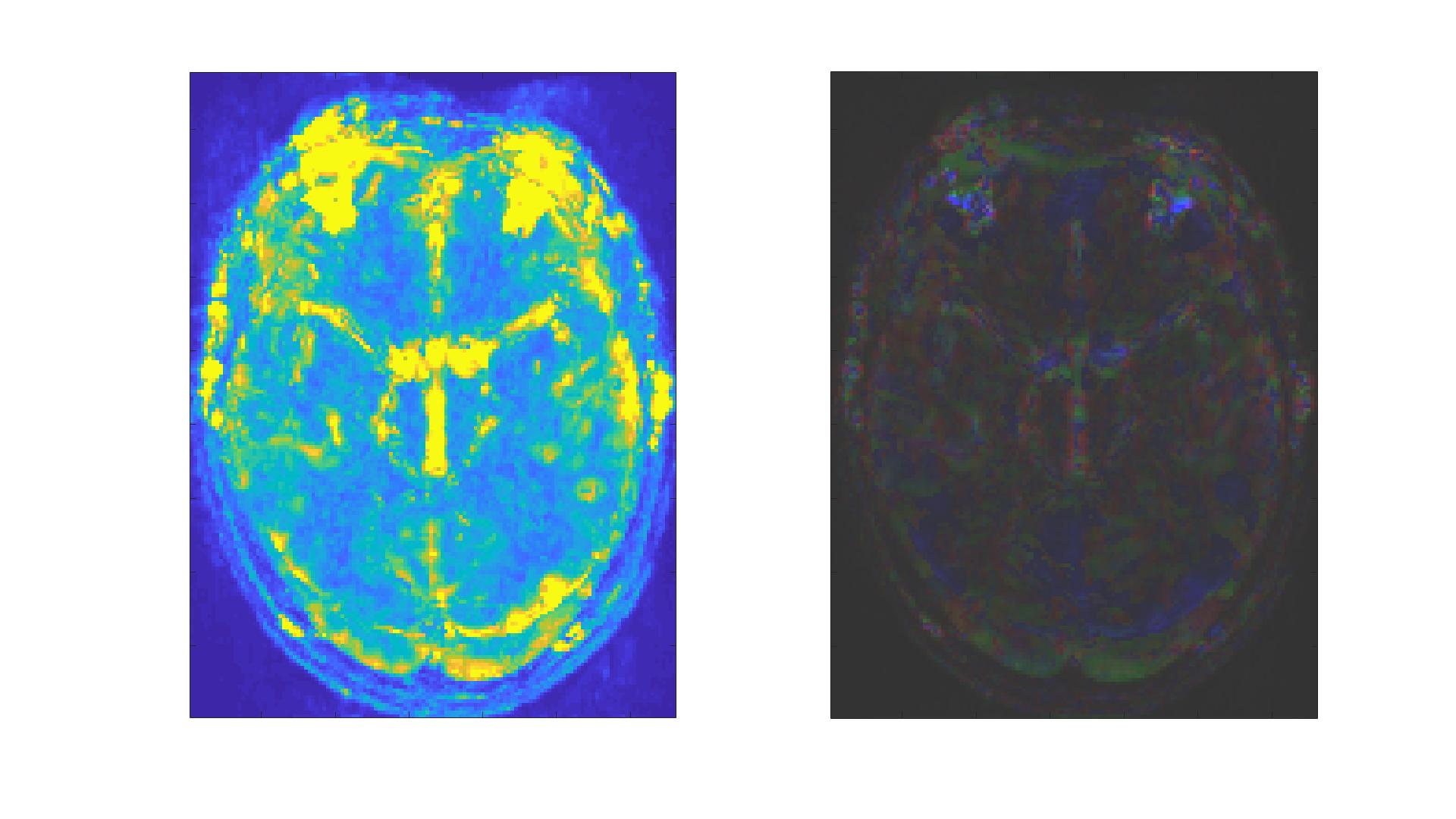}
\caption{Results of our proposed algorithm applied to subject 2 (upper two figures) and subject 16 (lower two figures) of the data set. Velocity norm MIPs (left) and colour direction MIPs (right).}
\label{fig_brain_N}
\end{figure}

Figure~\ref{fig_brain_N} shows the results of our CGNE algorithm terminated after $10$ iterations, in the same way as in \cite{Hubmer_Neubauer_Ramlau_Voss_2018}. The left images in the figure show maximum intensity projections (MIPs) over the $z$-axis of the norm of the reconstructed PWV vector field $v$. One can clearly see the location of some of the major blood vessels as well as (relative) information on the absolute velocity of the PWV. The right images show colour direction MIPs, created by assigning an RGB value to every pixel of the MIP by first identifying the voxel whose velocity norm entered the MIP at that pixel, and then taking the absolute values of the components $v_1$, $v_2$, and $v_3$ of the PWV at that voxel as the red, green, and blue values of the RGB triplet at that pixel, respectively. For example, a red pixel in the colour direction MIP indicates movement along the $x$-axis, a green pixel along the $y$-axis and a blue pixel along the $z$-axis. uniform scaling was applied to the resulting figures to enhance their colours. From this one can get some idea about the direction of the PWV. 

From a visual comparison with the results in \citeVH, one can see that our CGNE approach yields much better results than the one proposed in \citeV and leads to similar results as the algorithm proposed in \cite{Hubmer_Neubauer_Ramlau_Voss_2018}. However, it is much faster than the iterative algorithm from \cite{Hubmer_Neubauer_Ramlau_Voss_2018}, requiring only around $4$ instead of $14$ s per iteration on the same workstation (Intel(R) Xeon(R) CPU E5-1620 v4 @3.50GHz). Furthermore, it also scales much better with respect to the size of the input, as all calculations are explicit and no matrices need to be assembled or stored. This difference in speed makes it possible to process MRI data sets with a much higher spatio-temporal resolution.

\section{Conclusion and Outlook}\label{sect_conclusion}

In this paper we proposed a CGNE-based algorithm for solving the parameter estimation problem of MRAI. Based on a reformulation of the problem as a linear operator equation with both a noisy operator and a noisy right-hand side, this algorithm allows for an efficient numerical implementation which can also handle MRI data sets with a high spatio-temporal resolution, which will be available in the near future. Numerical experiments on experimental data show the competitiveness of our algorithm in comparison with existing reconstruction algorithms.

\section{Support}

S. Hubmer and R. Ramlau were (partly) funded by the Austrian Science Fund (FWF): F6805-N36, project 5 and W1214-N15, project DK8. 

H. Voss acknowledges support by the Nancy M. and Samuel C. Fleming Research Scholar Award in Intercampus Collaborations, Cornell University.

\bibliographystyle{plain}
{\footnotesize
\bibliography{mybib}
}

\end{document}